\title{Characterisations and Galois conjugacy of generalised Paley maps}

\author{Gareth A. Jones\\
School of Mathematics\\
University of Southampton\\
Southampton SO17  1BJ, U.K.\\
{\tt G.A.Jones@maths.soton.ac.uk}\\
Tel.~+44 (0)23 80593654\\
Fax.~+44 (0)23 80595147\\
}

\documentclass[12pt]{article}
\usepackage[latin1]{inputenc}
\usepackage{a4wide}
\usepackage{latexsym}
\usepackage{amsmath}
\usepackage{amssymb}
\usepackage{amsfonts}
\newtheorem{thm}{Theorem}[section]
\newtheorem{lemma}[thm]{Lemma}
\newtheorem{cor}[thm]{Corollary}
\newtheorem{prop}[thm]{Proposition}
\date{}
\begin{document} 

\maketitle

\begin{abstract}
A {\sl generalised Paley map\/} is a Cayley map for the additive group of a finite field $F$, with a subgroup $S=-S$ of the multiplicative group as generating set, cyclically ordered by powers of a generator of $S$. We characterise these as the orientably regular maps with orientation-preserving automorphism group acting primitively and faithfully on the vertices; allowing a non-faithful primitive action yields certain cyclic coverings of these maps. We determine the fields of definition and the orbits of the absolute Galois group ${\rm Gal}\,\overline{\bf Q}$ on these maps, and we show that if $(q-1)/(p-1)$ divides $|S|$, where $|F|=q=p^e$ with $p$ prime, then these maps are the only orientably regular embeddings of their underlying graphs; in particular this applies to the Paley graphs, where $|S|=(q-1)/2$ is even.
\end{abstract}



\noindent{\bf Running head:} Generalised Paley maps\\

\vfill\eject

\section{Introduction}
A map $\cal M$ on an oriented surface is {\sl orientably regular\/} if its orientation-preserving automorphism group ${\rm Aut}^+{\cal M}$ acts transitively on the arcs (directed edges) of $\cal M$. A standard problem in topological graph theory is that of classifying the orientably regular embeddings of a given class of arc-transitive graphs. This has been achieved for several classes, such as complete graphs $K_n$~\cite{JJ}, complete bipartite graphs $K_{n,n}$~\cite{J2}, cocktail party graphs $K_n\otimes K_2$ and dipoles $D_n$~\cite{NS1}, merged Johnson graphs $J(n,m)_I$~\cite{J1}, $n$-cubes $Q_n$ for $n$ odd~\cite{DKN}, and, according to a recent announcement, also $Q_n$ for $n$ even. Here we extend this to a class of arc-transitive graphs which includes the Paley graphs~\cite{Pal}. 

A Paley graph $P_q$ is the Cayley graph for the additive group of a field $F_q$ of order $q\equiv 1$ mod~$(4)$, where the generating set $S$ consists of the non-zero squares; choosing a generator $s$ of the cyclic group $S$ determines an orientably regular embedding of $P_q$ called a Paley map, described by White in~\cite[\S 16.8]{Whi}. In~\cite{LP}, Lim and Praeger have extended the definition of $P_q$ to generalised Paley graphs $P_q^{(n)}$ by allowing the generating set to be any subgroup $S=-S\cong C_n$ of the multiplicative group of a field $F_q$. In \S 2 we define analogous generalised Paley maps ${\cal M}_q(s)$, where $s$ generates $S$, and after considering some examples in \S 3 we show in Theorem~4.1 that if $n$ is divisible by $(q-1)/(p-1)$, where $q=p^e$ for a prime $p$, then these $\phi(n)/e$ maps are, up to isomorphism,  the only orientably regular embeddings of $P_q^{(n)}$. Theorem~5.1 characterises the maps ${\cal M}_q(s)$ as the only orientably regular maps $\cal M$ for which ${\rm Aut}^+{\cal M}$ acts primitively and faithfully on the vertices; removing the faithfulness condition allows central cyclic coverings of these maps, together with orientably regular dipole maps, classified by Nedela and \v Skoviera in~\cite{NS1}. The proofs of these results use basic properties of finite permutation groups, especially Frobenius groups.

Grothendieck's theory of {\it dessins d'enfants}~\cite{Gro, JS2} shows that maps on compact oriented surfaces correspond to algebraic curves defined over the field $\overline{\bf Q}$ of algebraic numbers. In \S 6 some recent results of Streit, Wolfart and the author~\cite{JSW} are used to determine the fields of definition of these generalised Paley maps and their orbits under the action of the absolute Galois group ${\rm Gal}\,\overline{\bf Q}$.

\section{Generalised Paley maps}

Let $F$ be the finite field $F_q$ of order $q=p^e$, where $p$ is prime, and let $n$ be a divisor of $q-1$, with $n$ even if $p>2$. The multiplicative group  $F^*=F\setminus\{0\}$ of $F$ is cyclic of order $q-1$, so it has a unique subgroup $S$ of order $n$. Our conditions imply that $-1\in S$, so $S=-S$; the relation $v\sim w$ on $F$ given by $v-w\in S$ is therefore symmetric, so it defines an undirected graph $P=P_q^{(n)}$ of order $q$ and valency $n$, with vertex set $F$ and edges $v\sim w$. Following Lim and Praeger~\cite{LP} we will call $P$ a {\sl generalised Paley graph}, since it generalises the Paley graph~\cite{Pal} which arises when $q\equiv 1$ mod~$(4)$ and $n=(q-1)/2$, so that $S$ is the group of squares in $F^*$. The following result is straightforward, so the proof is omitted (see~\cite[Theorem~2.2(1)]{LP}):

\begin{lemma}
The following are equivalent:
\begin{itemize}

\item the graph $P$ is connected;
\item $S$ generates the additive group $F$;
\item $S$ acts irreducibly on $F$, regarded as a vector space over its prime field $F_p$;
\item $e$ is the multiplicative order of $p$ {\rm mod}~$(n)$, the least $i\geq 1$ such that $p^i\equiv 1$ {\rm mod}~$(n)$.

\end{itemize}
\end{lemma}

When these conditions are satisfied, with $n$ even if $p>2$, the connected graph $P$ is the Cayley graph for the additive group $F$ with respect to its generating set $S$. In these circumstances, which we assume from now on, we will call $q=p^e$ and $n$ an {\sl admissible pair\/}.

If $s$ generates the cyclic group $S$ then the cyclic ordering $1, s, s^2, \ldots, s^{n-1}$ of $S$ gives a rotation $v+1, v+s, v+s^2, \ldots, v+s^{n-1}$ of the neighbours of each vertex $v$ in $P$. This defines a map ${\cal M}={\cal M}_q(s)$ which embeds $P$ in an oriented surface, so that this rotation of neighbours is induced by the local orientation around $v$. Since the cyclic ordering of $S$ is the restriction to $S$ of an automorphism $v\mapsto sv$ of the additive group $F$, this map is orientably regular (see~\cite[Theorem 16-27]{Whi}, for example). We will call $\cal M$ a {\sl generalised Paley map}, since if $P=P_q$ it is one of the Paley maps defined by White in~\cite[\S 16.8]{Whi}.

The group $A\Gamma L_1(q)$ consists of the transformations
\[v\mapsto av^{\gamma}+b\eqno(2.1)\]
of $F$ where $a, b\in F$, $a\neq 0$, and $\gamma$ is an element of the Galois group $\Gamma={\rm Gal}\,F\cong C_e$ of $F$, generated by the Frobenius automorphism $v\mapsto v^p$. The affine group $AGL_1(q)$ consists of those transformations $(2.1)$ with $\gamma=1$. Let $A=AGL_1^{(n)}(q) $ denote the subgroup of order $nq$ in $AGL_1(q)$ consisting of the affine transformations $(2.1)$ with $\gamma=1$ and $a\in S$. This acts faithfully as a group of orientation-preserving automorphisms of $\cal M$, and it permutes the arcs transitively, so ${\rm Aut}^+{\cal M}=A$. This group is a semidirect product of a normal subgroup $T\cong F\cong (C_p)^e$ consisting of the translations $v\mapsto v+b\;(b\in F)$, by a complement $A_0\cong S\cong C_n$ consisting of the automorphisms $v\mapsto av\; (a\in S)$ fixing $0$.

\medskip

\noindent{\bf Example 2.1. Maps of small valency.} The cases where $n\leq 2$ are slightly exceptional, so it is useful to deal with them briefly here. If $n=1$ our hypotheses imply that $q=2$ and $s=1$, so we obtain a single generalised Paley map ${\cal M}={\cal M}_2(1)$; this is the map $\{2,1\}$ in the notation of Coxeter and Moser~\cite{CM}, an embedding of the complete graph $K_2$ in the sphere, with ${\rm Aut}^+{\cal M}=AGL_1(2)\cong C_2$. If $n=2$ then $e=1$, $q$ is an odd prime $p$, and $s=-1$; the map ${\cal M}={\cal M}_p(-1)$ is the embedding $\{p,2\}$ of a cycle of length $p$ in the sphere, with two $p$-gonal faces, and ${\rm Aut}^+{\cal M}=AGL_1^{(2)}(p)$ is the dihedral group $D_p$ of order $2p$.

\begin{lemma}
Let $s$ and $s'$ be generators of $S$. Then ${\cal M}_q(s)\cong {\cal M}_q(s')$ if and only if $s$ and $s'$ are equivalent under the Galois group $\Gamma$ of $F_q$.
\end{lemma}

In order to prove this, and for further use later, we first summarise some basic general facts about orientably regular maps; for background, see~\cite{JS1}, for instance.

For any group $G$, the orientably regular maps $\cal M$ with ${\rm Aut}^+{\cal M}\cong G$ correspond to the generating pairs $x, y$ for $G$ such that $y$ has order $2$. Here $x$ is a rotation fixing a vertex $v$ of $\cal M$, sending each incident edge to the next incident edge according to the local orientation around $v$, while $y$ is a half-turn, reversing one of these incident edges, so that $z=(xy)^{-1}$ is a rotation preserving an incident face. We will call $x, y$ and $z$ the {\sl standard generators\/} of $G$. Conversely, given generators $x, y$ and $z$ of a group $G$ with $y^2=xyz=1$ one can construct a map $\cal M$, with arcs corresponding to the elements of $G$, and vertices, edges and faces corresponding to the cosets in $G$ of the cyclic subgroups generated by $x, y$ and $z$; this map has type $\{m,n\}$ where $m$ and $n$ are the orders of $z$ and $x$. Two such maps are isomorphic if and only if the corresponding sets of generators are equivalent under ${\rm Aut}\,G$.

Lemma~2.2 now follows immediately from the easily verified fact that if $n>1$ then ${\rm Aut}\,A$ can be identified with the group $A\Gamma L_1(q)$, acting by conjugation on its normal subgroup $A$.

Since there are $\phi(n)$ choices for a generator $s$ of $S$, permuted fixed-point-freely by $\Gamma$, it follows from Lemma~2.2 that there are, up to isomorphism, $\phi(n)/e$ maps ${\cal M}_q(s)$ for a given admissible pair $q$ and $n$. Generators $s$ and $s'$ of $S$ are equivalent under $\Gamma$ if and only if they have the same minimal polynomial over $F_p$; the $\phi(n)/e$ maps ${\cal M}_q(s)$ therefore correspond to the irreducible factors of the reduction mod~$(p)$ of the cyclotomic polynomial $\Phi_n(t)\in{\bf Z}[t]$ of the primitive $n$-th roots of $1$, all of degree $e$.

\begin{cor}
For a given admissible pair $q=p^e$ and $n$, the $\phi(n)/e$ generalised Paley maps ${\cal M}_q(s)$ are, up to isomorphism, the only orientably regular maps $\cal M$ of valency $n$ with ${\rm Aut}^+{\cal M}\cong AGL_1^{(n)}(q)$.
\end{cor}

\noindent{\sl Proof.} The case $n=1$ is trivial, so we may assume that $n>1$. Such maps $\cal M$ correspond to the orbits of ${\rm Aut}\,A=A\Gamma L_1(q)$ on generating pairs $x, y$ for $A=AGL_1^{(n)}(q)$ of orders $n$ and $2$. There are $\phi(n)q$ elements of order $n$ in $A$, namely the $\phi(n)$ generators of each of the $q$ mutually disjoint vertex stabilisers $A_v\cong C_n$ in $A$. If $p=2$ there are $q-1$ elements $y$ of order $2$, namely the non-identity elements of the translation group $T$; since $A_v$ acts irreducibly by conjugation on $T$, each of the $\phi(n)q(q-1)$ pairs $x, y$ generates $G$. If $p>2$ there are $q$ elements $y$ of order $2$, namely one in each vertex stabiliser $A_v$; in this case the maximality of $A_v$ in $A$ implies that $x$ and $y$ generate $A$ provided they are not in the same vertex stabiliser, so again there are $\phi(n)q(q-1)$ generating pairs $x, y$.

The automorphism group $A\Gamma L_1(q)$ of $A$ has order $eq(q-1)$ and it acts fixed-point-freely on these generating pairs, so it has $\phi(n)/e$ orbits on them. Thus there are $\phi(n)/e$ mutually non-isomorphic orientably regular maps $\cal M$ of valency $n$ with ${\rm Aut}^+{\cal M}\cong A$. Since this is the number of generalised Paley maps ${\cal M}_q(s)$, all of which have this valency and automorphism group, the result is proved. \hfill$\square$

\medskip

If $s$ and $s'$ are generators for $S$ then $s'=s^j$ for some $j$ coprime to $n$; thus for any admissible pair $q$ and $n$ the maps ${\cal M}_q(s)$ are all equivalent under Wilson's operations $H_j$, which raise local edge rotations to their $j$-th powers for $j$ coprime to the valency~\cite{Wil}. These operations form a group isomorphic to the group of units mod~$(n)$, and Lemma~2.2 shows that the stabiliser of any map ${\cal M}_q(s)$ (the {\sl exponent group\/} introduced by Nedela and \v Skoviera~\cite{NS2}) is the subgroup of order $e$ generated by $H_p$.

An orientably regular map $\cal M$ is said to be {\sl reflexible\/} if it has an orientation-reversing automorphism, so that it is isomorphic to its mirror image $H_{-1}({\cal M})$, and thus ${\rm Aut}^+{\cal M}$ has index $2$ in the full automorphism group ${\rm Aut}\,{\cal M}$; otherwise, $\cal M$ and $H_{-1}({\cal M})$ form a {\sl chiral pair}. Since $H_{-1}({\cal M}_q(s))={\cal M}_q(s^{-1})$, we see that ${\cal M}_q(s)$ is reflexible if and only if $s^{p^i}=s^{-1}$, or equivalently $p^i\equiv -1$ mod~$(n)$, for some $i=0, 1,\ldots, e-1$. If this holds then $p^{2i}\equiv 1$ mod~$(n)$, so either $n=2$ and $q=p>2$ (see Example~2.1), or $n>2$ and $e=2i$ is even with $p^{e/2}\equiv -1$ mod~$(n)$ (see the examples in \S3, for instance). 

If $\cal M$ is an orientably regular map then all its faces are $m$-gons for some $m$, and all its vertices have valency $n$ for some $n$; in the notation of~\cite{CM} we say that $\cal M$ has {\sl type\/} $\{m,n\}$. A {\sl Petrie polygon\/} in a map is a closed zig-zag path within the graph, turning alternately first left and first right at each successive vertex~\cite[\S 5.2]{CM}; in an orientably regular map these all have the same {\sl Petrie length\/}. 

\begin{lemma}
Let $\cal M$ be a generalised Paley map ${\cal M}_q(s)$ with $n>2$. If $n\equiv 0, 1$ or $3$ {\rm mod}~$(4)$ then $\cal M$ has type $\{n,n\}$ and genus $1+\frac{1}{4}q(n-4)$, whereas if  $n\equiv 2$ {\rm mod}~$(4)$ then $\cal M$ has type $\{n/2, n\}$ and genus $1+\frac{1}{4}q(n-6)$. The Petrie length is $2p$, where $q=p^e$.
\end{lemma}

[The cases $n\leq 2$ are exceptional; they were dealt with in Example~2.1.]

\medskip

\noindent{\sl Proof.} All vertices of $\cal M$ have valency $n=|S|$. Successive vertices around one particular face are $0,\; s,\; s-s^2,\; s-s^2+s^3,\; \ldots\;$,
so the face-valency $m$ is the least $j>0$ such that
\[s-s^2+s^3-\cdots-(-s)^j=0.\]
Now
\[(s^{-1}+1)(s-s^2+s^3-\cdots-(-s)^j)=1-(-s)^j,\]
with $s^{-1}+1\neq 0$ if $n>2$, so $m$ is the multiplicative order of $-s$. This is $n$ unless $2<n\equiv 2$ mod~$(4)$, in which case it is $n/2$. Since $\cal M$ has $q$ vertices, $nq/2$ edges and $nq/m$ faces, it has Euler characteristic
\[\chi=q\Bigl(1-\frac{n}{2}+\frac{n}{m}\Bigr),\]
which immediately gives the required genus $g=1-\frac{1}{2}\chi$. A typical Petrie polygon in $\cal M$ passes through vertices $1,\, 0,\, s,\, s-1,\, 2s-1,\, 2s-2,\, 3s-2,\, 3s-3,\, \ldots\,$ in that order, so the Petrie length is $2p$.\hfill$\square$

\section{Further examples}

\medskip

\noindent {\bf Example 3.1. Complete maps.} If $n=q-1$, so that $S=F^*$, then $P_q^{(n)}$ is the complete graph $K_q$, and the maps ${\cal M}_q(s)$ are the $\phi(n)/e$ orientably regular complete maps constructed by Biggs in~\cite{Big}, with $A= AGL_1(q)$. These have type $\{n/2, n\}$ if $q\equiv 3$ mod~$(4)$, and type $\{n, n\}$ otherwise. They are reflexible for $q\leq 4$, but for $q\geq 5$ they occur in chiral pairs. For instance, if $q=4$ there are two possible generators $s$ of $S=F_4^*\cong C_3$, conjugate under the Galois group of $F_4$, giving rise to a single orientably regular embedding $\cal M$ of $K_4$: this is the tetrahedral map $\{3,3\}$, a reflexible map on the sphere with ${\rm Aut}^+{\cal M}\cong A_4$ and ${\rm Aut}\,{\cal M}\cong S_4$. James and the author showed in~\cite{JJ} that these maps ${\cal M}_q(s)$ are the only orientably regular embeddings of any complete graphs; further details of these maps are given there, and also by White in~\cite[\S 16-4]{Whi}.

\medskip

\noindent {\bf Example 3.2. Paley maps.} If $q\equiv 1$ mod~$(4)$ and $n=(q-1)/2$, so that $S$ consists of the squares in $F^*$, then $P_q^{(n)}$ is the Paley graph $P_q$ introduced by Paley in~\cite{Pal}, and the maps ${\cal M}_q(s)$ are the $\phi(n)/e$ orientably regular Paley maps described by Biggs and White in~\cite[\S 5.7]{BW}, and by White in~\cite[\S 16-8]{Whi}. Only the unique Paley maps ${\cal M}_q(s)$ with $q=5$ and $q=9$ are reflexible.


If $5<q\equiv 5$ mod~$(8)$ then each Paley map ${\cal M}_q(s)$ has type $\{n/2,n\}$ and genus $(q^2-13q+8)/8$. For instance, if $q=13$ there are two such maps, namely the chiral pair of torus maps $\{3,6\}_{3,1}$ and $\{3,6\}_{1,3}$ described in~\cite[\S8.4]{CM}, corresponding to $s=4$ and $-3$ respectively. In the next case, with $q=29$, we find three chiral pairs of maps of type $\{7,14\}$ and genus $59$; these are denoted by C59.4, C59.5 and C59.6 in Conder's list of chiral maps~\cite{Con}.

If $q=5$ there is one Paley map ${\cal M}_5(-1)$, the embedding $\{5,2\}$ of the $5$-cycle $P_5$ on the sphere.

If $q\equiv 1$ mod~$(8)$ then each Paley map ${\cal M}_q(s)$ has type $\{n,n\}$ and genus $(q-1)(q-8)/8$. For instance, if $q=9$ there is one such map: we can take $F={\bf Z}_3[i]$ with $i^2=-1$, so $s=i=(1-i)^2$ generates $S$ (as does its Galois conjugate $i^3=-i$); the resulting map ${\cal M}_9(i)\cong{\cal M}_9(-i)$, illustrated in~\cite[\S 5.7]{BW}, is the reflexible torus map $\{4,4\}_{3,0}$ described in~\cite[\S8.3]{CM}. In the next case, with $q=17$, there are two chiral pairs of maps of type $\{8,8\}$ and genus $18$ (C18.1 and their duals in~\cite{Con}), and when $q=25$ there is one chiral pair of self-dual maps of type $\{12,12\}$ and genus $51$ (C51.16 in~\cite{Con}).

All other examples of Paley maps have genus $g>101$, so they do not appear in~\cite{Con}.

\medskip

\noindent{\bf Example 3.3. More maps of small valency.} The cases $n=1$ and $n=2$ were dealt with in \S 2. If $n=3$ then $q=4$ and we obtain the tetrahedral map $\{3,3\}$ mentioned in Example~3.1. If $n=4$ then $q=p$ or $p^2$ as $p\equiv 1$ or $3$ mod~$(4)$, and we respectively obtain the chiral torus maps $\{4,4\}_{a,b}$ with $p=a^2+b^2$, or the reflexible torus maps $\{4,4\}_{p.0}$. If $n=5$ then $q=16$ and we obtain the reflexible map of genus $5$ and type $\{5,5\}$ denoted by R5.9 in~\cite{Con}. If $n=6$ then $q=p$ or $p^2$ as $p\equiv 1$ or $5$ mod~$(6)$, and we respectively obtain the chiral torus maps $\{3,6\}_{a,b}$ with $p=a^2+ab+b^2$, or the reflexible torus maps $\{3,6\}_{p.0}$. If $n=7$ then $q=8$ and we obtain the Edmonds maps, a chiral pair of embeddings of $K_8$ of genus $7$ and type $\{7,7\}$ denoted by C7.2 in~\cite{Con}. If $n=8$ then $q=p$ or $p^2$ as $p\equiv 1$ mod~$(8)$ or $p\equiv 3, 5$ or $7$ mod~$(8)$, and we respectively obtain four or two maps of type $\{8,8\}$ and genus $1+q$; for $q=17$ we have the four Paley maps mentioned in Example~3.2, while for $q=25$ we have the chiral pair C26.1 in~\cite{Con}. If $n=9$ then $q=64$ and we obtain the reflexible map of type $\{9,9\}$ and genus $81$ denoted by R81.125 in~\cite{Con}.

\section{Characterisation of generalised Paley maps}

Let $q$ and $n$ be an admissible pair, so that $P=P_q^{(n)}$ is connected, and suppose that  $(q-1)/(p-1)$ divides $n$. (This includes the case where $q\equiv 1$ mod~$(4)$ and $n=(q-1)/2$, so that $P$ is the Paley graph $P_q$, and also all cases where $q=p$.) We aim to show that the only orientably regular embeddings of $P$ are the generalised Paley maps ${\cal M}_q(s)$ described in \S 2.

Lim and Praeger~\cite[Theorem 1.2(4)]{LP} have shown that  if $(q-1)/(p-1)$ divides $n$ then ${\rm Aut}\,P$ is the subgroup $H=A\Gamma L^{(n)}_1(q)$ of index $(q-1)/n$ in $A\Gamma L_1(q)$ consisting of those transformations $(2.1)$ with $a\in S$. This group $H$ is a semidirect product of an elementary abelian normal subgroup $T\cong(C_p)^e$, consisting of the translations $v\mapsto v+b\;(b\in F)$, by the stabiliser $H_0$ of the vertex $0$, consisting of the transformations $v\mapsto av^{\gamma}$ with $a\in S$ and $\gamma\in\Gamma$. Similarly $H_0$ is a semidirect product of a normal subgroup $D\cong S\cong C_n$, consisting of the transformations $v\mapsto av$ with $a\in S$, by a group $C\cong\Gamma\cong C_e$ of transformations $v\mapsto v^{\gamma}\;(\gamma\in\Gamma)$; it has a presentation
\[H_0=\langle c, d \mid c^e=d^n=1, d^c=d^p\rangle,\eqno(4.1)\]
where $c: v\mapsto v^p$ generates $C$ and $d$ generates $D$.

\begin{thm}
Let $q=p^e$ and $n$ be an admissible pair, with $n$ divisible by $(q-1)/(p-1)$. A map $\cal M$ is an orientably regular embedding of the generalised Paley graph $P^{(n)}_q$ if and only if $\cal M$ is isomorphic to a generalised Paley map ${\cal M}_q(s)$, where $s$ generates the subgroup $S$ of order  $n$ in $F_q^*$.
\end{thm}

\noindent{\sl Proof.} If $\cal M$ is any orientably regular embedding of a generalised Paley graph $P=P^{(n)}_q$ then the orientation-preserving automorphism group $G={\rm Aut}^+{\cal M}$ of $\cal M$ is a subgroup of $H={\rm Aut}\,P$ of order $nq$, with $G_0\cong C_n$ acting regularly on the set $S$ of neighbours of $0$. Since $G_0\leq H_0$ we look for elements $g$ of order $n$ in $H_0$ as possible generators $x$ for $G_0$.

\begin{lemma}
The only elements of order $n$ in $H_0$ are the generators of $D$.
\end{lemma}

\noindent{\sl Proof.} The presentation~$(4.1)$ shows that each element of $H_0$ has the form $g=d^ic^j$, with $i=0, 1, \ldots, n-1$ and $j=0, 1, \ldots, e-1$. Suppose that $g$ has order $n$. If $g\in D$ then $g$ is a generator of $D$, as required, so suppose that $g\not\in D$, that is, $j\neq 0$. The image of $g$ in $H_0/D\cong C_e$ has order $f=e/(j,e)=[j,e]/j$, which divides $n$. Replacing $g$ with a suitable primitive power, we may replace $j$ with $(j,e)$, that is, we may assume that $j$ divides $e$, so $fj=e$. Induction on $k$ shows that
\[g^k=c^{kj}d^{i(p^{kj}+p^{(k-1)j}+\cdots+p^{2j}+p^j)}\]
for each $k\geq 0$. Since $c^{fj}=1$ we have
\[g^f=d^{i(p^{fj}+p^{(f-1)j}+\cdots+p^{2j}+p^j)}.\]
Since $g$ has order $n$, which is divisible by $f$, it follows that $g^f$ has order $n/f$. However, as a power of $d$ its order is equal to $n/(i(p^{fj}+p^{(f-1)j}+\cdots+p^{2j}+p^j),n)$, so
\[(i(p^{fj}+p^{(f-1)j}+\cdots+p^{2j}+p^j),n)=f.\eqno(4.2)\]
Now $p^{(f-1)j}+\cdots+p^j+1$ divides $p^{fj}+p^{(f-1)j}+\cdots+p^{2j}+p^j$, and it also divides
\[\frac{(p^j-1)}{(p-1)}(p^{(f-1)j}+\cdots+p^j+1)=\frac{p^{fj}-1}{p-1}=\frac{p^e-1}{p-1}=\frac{q-1}{p-1},\]
which by our hypothesis divides $n$.
Clearly $p^{(f-1)j}+\cdots+p^j+1>f$, giving a contradiction to $(4.2)$.
Thus the only elements of order $n$ in $H_0$ are the generators of $D$.\hfill$\square$

\medskip

Recall that $A=AGL_1^{(n)}(q)$, the group of transformations $(2.1)$ with $\gamma=1$ and  $a\in S$.

\begin{cor}
$G=A$.
\end{cor}

\noindent{\sl Proof.} The case $n=1$ is trivial (see Example~2.1), so we may assume that $n>1$. Lemma~4.2 implies that $G_0=D$. Since $\cal M$ is orientably regular, $G$ acts transitively on the vertices $v$ of $\cal M$, so their stabilisers $G_v$  in $G$ are the conjugates of $D$ in $G$. Now $AGL_1(q)$ is a normal subgroup of $A\Gamma L_1(q)$, so $A=H\cap AGL_1(q)$ is a normal subgroup of $H$. We have seen that $G_0=D\leq A$, and we have $G\leq H$, so $G_v\leq A$ for every vertex $v$.

For any orientably regular map, the vertex stabiisers $G_v$ generate a subgroup of index at most $2$ in $G$, namely the normal closure of $x$. If this index is $2$ then the embedded graph is bipartite, since the map covers the $2$-vertex map $\{2,1\}$. This is not the case here since ${\rm Aut}\,P$ acts primitively on the vertices for $n>1$. Thus $G$ is generated by the stabilisers $G_v$, so $G\leq A$. Since $|G|=nq=|A|$ it follows that $G=A$.\hfill$\square$


\medskip

The proof of Theorem~4.1 now follows immediately from Corollaries~2.3 and 4.3. \hfill$\square$

\begin{cor}
The only orientably regular embeddings of the Paley graphs $P_q$ are the $\phi(\frac{1}{2}(q-1))/e$ Paley maps ${\cal M}_q(s)$, where $q=p^e$ with $p$ prime, and $s$ generates the group of squares in $F_q^*$.\hfill$\square$
\end{cor}

\noindent{\bf Remarks. 1.} In order to show that  ${\rm Aut}\,P=H$, this proof of Theorem~4.1 relies on Theorem~1.2(4) of~\cite{LP}, which in turn relies on the classification of finite simple groups. However there is a more direct proof of Corollary~4.4, instead using a result of Carlitz~\cite{Car} (see also~\cite{Muz}) which gives an elementary proof that ${\rm Aut}\,P=H$ for the Paley graphs $P=P_q$. Similarly in the cases $q=p$ and $p^2$ of Theorem~4.1 one can use classical results of Burnside~\cite{Bur1} (see also~\cite[\S 251]{Bur2}) and Wielandt~\cite[\S 16]{Wie} on primitive permutation groups of these degrees.

\smallskip

\noindent{\bf 2.} As pointed out by Lim and Praeger in~\cite{LP}, ${\rm Aut}\,P$ can be much larger than $A\Gamma L_1(q)$ if $n$ is not divisible by $(q-1)/(p-1)$. For instance, if $q=p^2$ and $n=2(p-1)$ then $P$ is a Hamming graph $H(2,p)$, and ${\rm Aut}\,P$ is the wreath product $S_p\wr S_2$ of order $2(p!)^2$. In this particular case it is not hard to show that the only orientably regular embeddings are again the generalised Paley maps ${\cal M}_q(s)$, but the general problem of classifying orientably regular embeddings remains open, for Hamming graphs and for generalised Paley graphs.



\section{Vertex-primitive maps}

Here we consider the orientably regular maps $\cal M$ for which ${\rm Aut}^+{\cal M}$ acts primitively on the vertices; of course, this includes all cases where the number of vertices is prime. First we classify the maps for which the action on vertices is also faithful.

\begin{thm}
Let $\cal M$ be an orientably regular map on a compact surface. Then the orientation-preserving automorphism group ${\rm Aut}^+{\cal M}$ acts primitively and faithfully on the vertices of $\cal M$ if and only if $\cal M$ is isomorphic to a generalised Paley map ${\cal M}_q(s)$.
\end{thm}

\noindent{\sl Proof.} Each map ${\cal M}_q(s)$ has these properties since its automorphism group $AGL_1^{(n)}(q)$ acts faithfully on $F$, with the stabiliser of $0$ acting as an irreducible linear group.

Conversely, if the group $G={\rm Aut}^+{\cal M}$ acts primitively on the vertex set $V$ of a map $\cal M$, then each vertex stabiliser $G_v$ is a maximal subgroup of $G$. If $G_v=1$ then $G$ is cyclic of prime order, and since $G$ contains an involution we must have $G\cong C_2$, so ${\cal M}\cong {\cal M}_2(1)$, the planar embedding $\{2,1\}$ of $K_2$. We may therefore assume that $G_v\neq 1$. If $v$ and $w$ are distinct vertices, then $G_{vw}=1$: for if $g\in G_{vw}=G_v\cap G_w$ then since the vertex stabilisers are abelian, distinct and maximal, the centraliser $C_G(g)$ of $g$ contains $\langle G_v, G_w\rangle=G$, so $g$ is in the centre of $G$; now $g$ fixes at least one vertex, and the vertex stabilisers are all conjugate in $G$, so $g$ must fix every vertex, giving $g=1$.

This shows that $G$ acts on $V$ as a Frobenius group, so it has a Frobenius kernel $K$, a regular normal subgroup consisting of the fixed-point-free elements together with $1$ (see~\cite[\S V.8]{Hup} for properties of Frobenius groups). We may identify $V$ with $K$, acting by multiplication on itself. The stabiliser of the identity vertex is a complement for $K$ in $G$, and its action on $V$ coincides with its action by conjugation on $K$. Since $G$ acts primitively, no proper subgroup of $K$ can be normal in $G$, so $K$ is characteristically simple, that is, a direct product of isomorphic simple groups. By a theorem of Thompson~\cite{Tho}, finite Frobenius kernels are nilpotent, so $K$ is an elementary abelian $p$-group for some prime $p$.

We can therefore regard $V$ as a vector space of some dimension $e$ over the field $F_p$ of order $p$. The subgroup $G_0$ fixing the vertex $0$ acts as a group of linear transformations of $V$, and the primitivity of $G$ implies that $V$ is an irreducible $G_0$-module. Since $G_0$ is cyclic we can therefore identify $V$ with the field $F_q$ of order $q=p^e$ so that $G_0$ acts by multiplication as a subgroup $S$ of $F_q^*$ (see~\cite[Satz II.3.10]{Hup}, for instance). All orbits of $S$ on $V\setminus\{0\}$ have length $n=|S|$, and since $\cal M$ is orientably regular the neighbours of $0$ form such an orbit, so $\cal M$ has valency $n$. Since $|G|=nq$ must be even, $n$ is even if $p>2$. The irreducibility of $V$ implies that $S$ is not contained in any proper subfield of $F_q$, so $q$ and $n$ form an admissible pair. Since $G=VG_0=AGL_1^{(n)}(q)$ it follows from Corollary~2.3 that $\cal M$ is isomorphic to a generalised Paley map ${\cal M}_q(s)$ for some $s$. \hfill $\square$

\medskip

The condition that the group $G={\rm Aut}^+{\cal M}$ should act faithfully on the $q$ vertices of $\cal M$ is not particularly restrictive. (Li and \v Sir\'a\v n discuss non-faithful actions on vertices, edges and faces for more general maps in~\cite{LS}.) If this action is primitive but not faithful, then the kernel of the action (the intersection of the vertex stabilisers) is a normal subgroup $N\cong C_k$ of $G$ where $k$ is the number of edges joining each pair of adjacent vertices.  By Theorem~5.1 the quotient map $\overline{\cal M}={\cal M}/N$ is a generalised Paley map ${\cal M}_q(s)$, with orientation-preserving automorphism group $\overline G=G/N\cong AGL_1^{(n)}(q)$ for some admissible pair $q$ and $n$, and $\cal M$ and $G$ are $k$-fold cyclic coverings of $\overline{\cal M}$ and $\overline G$.

Let $\overline{\cal M}={\cal M}_q(s)$ have type $\{m,n\}$, and let $x, y$ and $z$ be the standard generators of $G$. Then $N=\langle x^n\rangle\cong C_k$ and $z^m\in N$, so $z^m=x^{in}$ for some $i$. The covering ${\cal M}\to\overline{\cal M}$ is branched over the vertices of $\overline{\cal M}$, and also over the faces if $i\not\equiv 0$ mod~$(k)$. There are $q$ vertices and $knq/2$ edges in $\cal M$, and since $z$ has order $km/(i,k)$ there are $(i,k)nq/m$ faces. Thus $\cal M$ has type $\{km/(i,k), kn\}$ and genus
\[g=1+\frac{q}{4m}\bigl(kmn-2m-2(i,k)n\bigr).\]

\medskip

\noindent{\bf Example 5.1.} If an integer $u$ satisfies $u^2\equiv 1$ mod~$(k)$ then the group
\[ G=\langle x, y \mid x^k=y^2=1, x^y=x^u\rangle\]
of order $2k$ has the form ${\rm Aut}^+{\cal M}$ where $\cal M$ is the reflexible {\sl dipole map\/} ${\cal D}_k(u)$, with two vertices of valency $k$. (Nedela and \v Skoviera showed in~\cite{NS1} that these are the only orientably regular embeddings of dipole graphs.) There is a normal subgroup $N=\langle x\rangle\cong C_k$ in $G$ with $\overline G=G/N\cong  AGL_1(2)\cong C_2$, so $\cal M$ is a $k$-sheeted cyclic covering of the spherical embedding $\overline{\cal M}={\cal M}/N={\cal M}_2(1)=\{2,1\}$ of $K_2$. If $k>1$ there is a non-faithful action of $G$ on the two vertices of $\cal M$, with kernel $N$. The number of possible values $u\in{\bf Z}_k$, and hence of maps $\cal M$, is $2^{\mu+\nu}$, where $\nu$ is the number of distinct odd primes dividing $k$, and $\mu=2, 1$ or $0$ as $k\equiv 0, 4$ or otherwise mod~$(8)$. We have $z^2=(xy)^{-2}=x^{-u-1}$, so $i\equiv-u-1$ mod~$(k)$. The covering ${\cal M}\to\overline{\cal M}$ is branched over the two vertices of $\overline{\cal M}$, and also over its single face if $u\not\equiv -1$ mod~$(k)$; $\cal M$ has type $\{m,k\}$ where $m=2k/(u+1,k)$ is the order of $z$, and its genus is $(k-(u+1,k))/2$.  For instance, if $k=8$ then $u\equiv 1, 3, 5$ or $7$ mod~$(8)$, and $\cal M$ is respectively R3.11, R2.3, R3.10 in~\cite{Con}, or the spherical map $\{2,8\}$.

\begin{prop}
If $G={\rm Aut}^+{\cal M}$ acts primitively on the vertices of $\cal M$, and the kernel $N$ is not contained in the centre of $G$, then ${\cal M}\cong {\cal D}_k(u)$ for some $u\not\equiv 1$ {\rm mod}~$(k)$.
\end{prop}

\noindent{\sl Proof.}  Let $x, y$ and $z$ be the standard generators of $G$. Since $N$ is contained in all the vertex stabilisers, and these are the conjugates of $\langle x\rangle$, its centraliser $C=C_G(N)$ contains the normal closure of $x$ in $G$, so $C$ has index at most $2$ in $G$. By our hypothesis $C\neq G$, so the index is $2$ and hence the underlying graph of $\cal M$ is bipartite, as in the proof of Corollary~4.3. Since the vertices are permuted primitively there must be just two of them, so the graph is a dipole and hence ${\cal M}\cong{\cal D}_k(u)$ for some $u$ with $u^2\equiv 1$ mod~$(k)$ by~\cite{NS1}. Since $G$ is nonabelian we have $u\not\equiv 1$ mod~$(k)$. \hfill$\square$

\medskip

This result directs attention towards the cases such as ${\cal D}_k(1)$ where $N$ is in the centre of $G$. The monodromy permutation of the sheets of the covering is then $x^n$ at each vertex, and $z^{m}=x^{in}$ at each face. A closed path in $\overline{\cal M}$ going once round each of these branch-points is homologically trivial, so the product of these monodromy permutations must be the identity. There are $q$ vertices and $f=nq/m$ faces in $\overline{\cal M}$, so this implies that
\[q+fi\equiv 0\; {\rm mod}\,(k).\eqno(5.1)\]
This has $(f,k)$ solutions $i\in{\bf Z}_k$ if $(f,k)$ divides $q$, and none otherwise, giving an upper bound on the number of $k$-sheeted central coverings $\cal M$ for a given pair ${\cal M}_q(s)$ and $k$.

The following construction gives at least one example of such a map $\cal M$ for each pair ${\cal M}_q(s)$ and $k$, provided $k$ is odd when $p>2$. We can regard $\overline G$ as a semidirect product of the normal translation group $T\cong F$ by $\overline G_0\cong C_n$, with the standard generator of the complement $\overline G_0$ acting by conjugation on $T$ as a specific automorphism of order $n$. We form a semidirect product $G$ of $T$ by $G_0\cong C_{nk}$, with a generator $x$ of $G_0$ inducing the same automorphism of $T$, so that $x^n$ generates a central subgroup $N\cong C_k$ in $G$ with $G/N\cong \overline G$. For $G$ to be associated with an orientably regular map $\cal M$ that covers $\overline{\cal M}$ we need an involution $y\in G$ which projects onto the corresponding standard generator of $ \overline G$. If $p=2$ then we can use the same element $y\in T$ for $G$ as for $\overline G$, but if $p>2$ (so that $n$ is even) then we need $y=x^{jn+n/2}t$ for some $j\in{\bf Z}_k$ and non-identity $t\in T$. We have
\[y^2=x^{(2j+1)n}x^{-n/2}tx^{n/2}t=x^{(2j+1)n}\]
since $x^{n/2}$ inverts $T$, so we need $2j\equiv -1$ mod~$(k)$. Thus if $p>2$ then $k$ must be odd and $j\equiv  (k-1)/2$ mod~$(k)$, giving $y=x^{kn/2}t$. In either case, $G=\langle x,y\rangle$ since $G_0$ is maximal.

\medskip

\noindent{\bf Example 5.2.} Let $\overline{\cal M}={\cal M}_7(s)$ with $s=3$ or $5$, one of a chiral pair of torus embeddings of $K_7$, of type $\{3,6\}$. Here
\[\overline G=\langle x, t \mid x^6=t^7=1,\, t^x=t^s\rangle,\]
a semidirect product of $T=\langle t\rangle\cong C_7$ by $\overline G_0=\langle x\rangle\cong C_6$, with $y=x^3t$. We can define
\[G=\langle x, t \mid x^{6k}=t^7=1,\, t^x=t^s\rangle,\]
a semidirect product of $\langle t\rangle\cong C_7$ by $\langle x\rangle\cong C_{6k}$. This has centre $N=\langle x^6\rangle\cong C_k$, with $G/N\cong\overline G$. By the preceding argument we can take $y=x^{3k}t$ with $k$ odd. Then $xy=x^{3k+1}t$ giving $z^3=x^{-3(k+1)}$, so $z$ has order $m=3k$. Thus for each odd $k$ there is a chiral pair of orientably regular maps $\cal M$ of type $\{3k, 6k\}$ and genus $(21k-19)/2$ with ${\rm Aut}^+{\cal M}=G$ and ${\cal M}/N\cong\overline{\cal M}$; each has seven vertices, permuted primitively by $G$, with kernel $N$. For $k=3, 5, 7$ and $9$ these are the chiral pairs C22.6, C43.10, C64.20 and C85.14 in~\cite{Con}. There are no such central coverings when $k$ is even, since $(5.1)$ gives $7+14i\equiv 0$ mod~$(k)$.

\medskip

\noindent{\bf Example 5.3.} Let $\overline{\cal M}$ be the reflexible spherical map ${\cal M}_p(-1)=\{p,2\}$, with $\overline G=AGL_1^{(2)}(p)\cong D_p$. Thus $m=q=p>2$ and $n=2$. Since $f=2$, $(5.1)$ gives $p+2i\equiv 0$ mod~$(k)$. If $k$ is odd there is a unique solution $i\equiv (k-p)/2$ mod~$(k)$, giving a map $\cal M$ with
\[G=\langle x, t \mid x^{2k}=t^p=1,\, t^x=t^{-1}\rangle\cong D_p\times C_k\]
where $y=x^kt$ and $z=(x^{k+1}t)^{-1}$, and $N=\langle x^2\rangle$; this map has type $\{k,2k\}$ or $\{kp,2k\}$ as $p$ divides $k$ or not, and its genus is $1+\frac{1}{2}(k-3)p$ or $\frac{1}{2}(k-1)p$ respectively. If $k$ is even there is no solution of $(5.1)$, and hence no map. For instance, let $k=3$: for $p=3$ we have the torus map $\{3,6\}_{1,1}$ with six triangular faces, and for $p>3$ a map of type $\{3p,6\}$ and genus $p$; for $p=5$ and $7$ these are the duals of R5.11 and R7.8 in~\cite{Con}.

\medskip

These examples have $p>2$, so the construction fails for even $k$, and $(5.1)$ shows that no such coverings can exist. The next example, with $p=2$, yields a covering for every $k$:

\medskip

\noindent{\bf Example 5.4.} Let $\overline{\cal M}={\cal M}_2(1)$, the reflexible spherical embedding $\{2,1\}$ of $K_2$, with $\overline G=AGL_1(2)\cong C_2$. Here $q=m=2$, $n=1$ and $f=1$, so $(5.1)$ gives $2+i\equiv 0$ mod~$(k)$; this has a unique solution $i\equiv -2$ mod~$(k)$ for each $k\geq 1$, corresponding to a unique map $\cal M$ with $G=\langle x, y \mid x^k=y^2=[x,y]=1\rangle\cong C_2\times C_k$. This is the reflexible dipole map ${\cal D}_k(1)$ in Example~5.1. If $k$ is even it has type $\{k,k\}$ and genus $(k-2)/2$; if $k$ is odd it has type $\{2k,k\}$ and genus $(k-1)/2$. For $k=2$ we have the spherical map $\{2,2\}$ of two digons. For $k=3$ we have the torus map $\{6,3\}_{1,0}$ with a single hexagonal face. For $k=4$ we have the torus map $\{4,4\}_{1,1}$ with two square faces. For $k=5$ and $6$ we have the maps $\{10,5\}_2$ and $\{6,6\}_2$ of genus $2$ in~\cite[Table 9]{CM}; in~\cite{Con}, the first is the dual of the map R2.4, and the second is R2.5. For $k=7$ we have the dual of R3.9, and for $k=8$ we have R3.11.

\medskip

The last three examples may give the impression that a given pair $\overline{\cal M}$ and $k$ yields at most one central covering $\cal M$. The following example shows that this is not  always true:

\medskip

\noindent{\bf Example 5.5.} Let $\overline{\cal M}={\cal M}_4(s)$, the spherical embedding $\{3,3\}$ of a tetrahedron, with $\overline G=AGL_1(4)\cong A_4$. Here $q=f=4$, so $(5.1)$ gives $4(1+i)\equiv 0$ mod~$(k)$, with $(4,k)$ solutions $i\in{\bf Z}_k$. For instance, if $k=4$ there are four maps $\cal M$, namely R3.3 of type $\{3,12\}$ and genus $3$, R7.7 of type $\{6,12\}$ and genus $7$, and R9.26 and R9.27 of type $\{12,12\}$ and  genus $9$; if $k=8$ there are four maps R21.32 -- R21.35, all of type $\{24,24\}$ and genus $21$.

\section{Galois actions}

According to Grothendieck's theory of {\it dessins d'enfants}~\cite{Gro, JS2}, a map $\cal M$ on a compact oriented surface corresponds naturally to a {\sl Bely\u{\i} pair\/} $(X,\beta)$, where $X$ is a nonsingular projective algebraic curve over $\bf C$, and $\beta$ is a rational function from $X$ to the complex projective line (or Riemann sphere) ${\bf P}^1({\bf C})={\bf C}\cup\{\infty\}$, unramified outside $\{0, 1, \infty\}$. One can regard $X$ as a Riemann surface underlying $\cal M$, with the inverse image under $\beta$ of the unit interval providing the embedded graph. For instance, the dual of the map $\cal M$ in Example~5.3, with $k=p$, is the standard embedding of the complete bipartite graph $K_{p,p}$ described by Biggs and White in~\cite[\S 5.6.7]{BW}; here $X$ is the Fermat curve $x^p+y^p=z^p$, with $\beta([x,y,z])=(x/z)^p$, $\beta^{-1}(0)$ and $\beta^{-1}(1)$ giving the black and white vertices, and $\beta^{-1}([0,1])$ the edges~\cite{J2, JS2}. Bely\u{\i}'s Theorem~\cite{Bel} asserts that $X$ and $\beta$ are defined (by polynomials and rational functions) over the field $\overline{\bf Q}$ of algebraic numbers, and as Grothendieck observed, the action of the {\sl absolute Galois group\/} ${\bf \Gamma}={\rm Gal}\,\overline{\bf Q}$ on their coefficients induces a faithful action of this group on the associated maps $\cal M$. Finding explicit fields of definition and Galois orbits is an important but usually difficult problem. The following result generalises some examples given by Streit, Wolfart and the author in~\cite{JSW}:

\begin{thm} For any admissible pair $q=p^e$ and $n$, the $\phi(n)/e$ generalised Paley maps ${\cal M}_q(s)$ form an orbit under $\bf \Gamma$, and the corresponding Bely\u{\i} pairs are defined over the splitting field of $p$ in the cyclotomic field ${\bf Q}(\zeta_n)$,  where $\zeta_n=\exp(2\pi i/n)$.
\end{thm}

[This field is the unique subfield of  ${\bf Q}(\zeta_n)$ of degree $\phi(n)/e$ over $\bf Q$.]

\medskip

\noindent{\sl Proof.} As shown by Streit and the author in~\cite{JSt}, the automorphism group ${\rm Aut}^+{\cal M}$ of a map $\cal M$ and various parameters such as its vertex-valencies are invariant under the action of $\bf \Gamma$. Corollary~2.3 shows that the set of generalised Paley maps ${\cal M}_q(s)$ for a given admissible pair $q, n$ is characterised by their common automorphism group and valency, so this set is $\bf \Gamma$-invariant. As noted in \S 2, these maps are all equivalent under Wilson's operations $H_j$. It therefore follows immediately from Theorem~2 of~\cite{JSW} that they form an orbit under $\bf \Gamma$, and that the corresponding Bely\u{\i} pairs are defined over the splitting field of $p$ in ${\bf Q}(\zeta_n)$. \hfill$\square$

\medskip

\noindent{\bf Example 6.1.} If we take $q=29$ and $n=14$ then the resulting $\phi(14)/1=6$ Paley maps ${\cal M}_{29}(s)$ of genus $59$ form an orbit under $\bf \Gamma$, and the corresponding Bely\u{\i} pairs are defined over ${\bf Q}(\zeta_{14})={\bf Q}(\zeta_7)$, with $\bf \Gamma$ permuting them as the Galois group $C_6$ of this field.






\end{document}